\documentclass[journal]{IEEEtran}

\usepackage{dcolumn}
\usepackage{bm}
\usepackage[noadjust]{cite}
\usepackage{amsmath}
\usepackage{array}
\usepackage{xcolor}

\usepackage{booktabs}
\usepackage{amssymb}
\usepackage{neuralnetwork}
\usepackage{tikz}
\usepackage{mwe}
\usepackage{graphbox}
\usepackage{graphicx,pdflscape}
\usepackage{rotating}
\usepackage{hyperref}
\usepackage{cleveref}
\usepackage{adjustbox}
\usetikzlibrary{matrix,chains,positioning,decorations.pathreplacing,arrows}
\usetikzlibrary{positioning}
\usepackage{algorithm}
\usepackage{algpseudocode}
\usepackage{xcolor}

\newcommand*\Let[2]{\State #1 $\gets$ #2}

\newcommand{\rev}[1]{{#1}}
\newcommand{\revtwo}[1]{{#1}}

\usepackage{xpatch}
\makeatletter
\xpatchcmd{\linklayers}{\nn@lastnode}{\lastnode}{}{}
\xpatchcmd{\linklayers}{\nn@thisnode}{\thisnode}{}{}
\makeatother

\begin{document}


\title{An efficient Quasi-Newton method for nonlinear inverse problems via learned singular values}

\author{Danny Smyl, Tyler N. Tallman, Dong Liu,~\IEEEmembership{Senior Member,~IEEE}, Andreas Hauptmann,~\IEEEmembership{Member,~IEEE,}
\thanks{This work was supported by Academy of Finland Proj. 336796, 334817, and the CMIC-EPSRC platform grant (EP/M020533/1), DL was supported by National Natural Science Foundation of China (Grant No. 61871356).}%
\thanks{D. Smyl is with the Department of Civil and Structural Engineering, University of Sheffield, UK }%
\thanks{T.N. Tallman is with the School of Aeronautics and Astronautics, Purdue University, West Lafayette, IN, USA}%
\thanks{D. Liu is with CAS Key Laboratory of Microscale Magnetic Resonance and Department of Modern Physics, University of Science and Technology of China, Hefei, 230026, China}
\thanks{A. Hauptmann is with the Research Unit of Mathematical Sciences; University of Oulu, Oulu, Finland and with the Department of Computer Science; University College London, London, United Kingdom.}%
}%

\maketitle

             

\begin{abstract}
Solving complex optimization problems in engineering and the physical sciences requires repetitive computation of multi-dimensional function derivatives\rev{, which commonly} require computationally-demanding numerical differentiation such as perturbation techniques. 
In particular, \rev{Gauss-Newton methods are used for} nonlinear inverse problems that require iterative updates to be computed from the Jacobian \rev{and allow for flexible incorporation of prior knowledge.}
Computationally more efficient alternatives are Quasi-Newton methods, where the repeated computation of the Jacobian is replaced by an approximate update\rev{, but unfortunately are often too restrictive for highly ill-posed problems.}
\rev{To overcome this limitation,} we present a highly efficient data-driven Quasi-Newton method applicable to nonlinear inverse problems, by using the singular value decomposition and learning a mapping from model outputs to the singular values to compute the updated Jacobian. \rev{Enabling} time-critical applications and allowing for flexible incorporation of prior knowledge necessary to solve ill-posed problems. We present results for the highly non-linear inverse problem of electrical impedance tomography with experimental data. 

\end{abstract}

%



\section{Introduction} 
Historically, derivative-based optimization is at the heart of, and required for solving, complex systems in  engineering, applied sciences, and applied mathematics \cite{belegundu2019optimization,biegler2014multi}.
The computation of suitable update equations can be particularly problematic, when the system equations are computationally expensive to evaluate. This is especially a problem, when solving (large) nonlinear inverse problems of the form
\begin{equation}\label{eqn:nonLinProb}
    g=\mathcal{A}(f)+\delta g,
\end{equation}
with a nonlinear model $\mathcal{A}$, typically given by a finite element solver. Given data $g$ with possible noise $\delta g$, we can then compute a solution $f$ minimizing a suitable penalty functional
\begin{equation}\label{eqn:leastSquaresReg}
    \mathcal{L}(f) = \|g - \mathcal{A}(f)\|_2^2 + R(f).
\end{equation}
The first data-fidelty term ensures closeness of solutions to the measurements and the second term incorporates necessary regularization to counteract the ill-posedness and stabilize the reconstruction procedure as well as incorporate prior information.  
For nonlinear problems, second-order methods are typically employed to compute minimizers of \cref{eqn:leastSquaresReg} due to fast (quadratic) convergence. \rev{However}, given potentially expensive evaluations of the forward model $\mathcal{A}$ the computation of second order derivatives is undesirable. Thus, approximations such as the Gauss-Newton method are used, where only first order derivatives are required in each iteration, collected in the Jacobian
$\textbf{J} = \frac{\partial \mathcal{A}(f) }{\partial f}$.
We can then form an approximate Hessian $\textbf{H} = \textbf{J}^T\textbf{J}$, characterizing the GN approach.
\rev{Yet,} if no efficient/accurate analytical or semi-analytical algorithms are available, computing $\textbf{J}$ is often problematic, \rev{requiring} either (a) numerical differentiation or (b)  updating approaches.

Both, (a) and (b), have advantages and disadvantages.
Numerical differentiation approaches, namely the perturbation method, can indeed yield reliable derivative approximations.
However, perturbation methods require repetitive function evaluations which can be extremely demanding.
\rev{Therefore, when} semi-analytical gradients are not available, conventional numerical differentiation approaches may be infeasible \rev{in solving high-dimensional nonlinear inverse problems} - if not intractable on accessible computing resources \cite{saibaba2013flexible,smyl2019less}.
To address this challenge, a number of researchers have proposed methods for updating gradient terms, rather than, e.g., computing finite differences at each iteration; these approaches are referred to as Quasi-Newton (QN) methods.
Possibly the most classical QN derivative approximator utilizes Broyden's method proposed in 1965 \cite{broyden1965class}, taking advantage of the secant method.
Since then, a number of successful updating methods have been proposed including the ever popular Broyden–Fletcher–Goldfarb–Shanno (BFGS) algorithm \cite{broyden1970convergence,fletcher1970new,goldfarb1970family,shanno1970conditioning}, the limited memory BFGS algorithm \cite{liu1989limited}, and a suite of other QN-family derivative approximators \cite{lewis2013nonsmooth}. Nevertheless, incorporation of prior information to regularize and stabilize the reconstruction for ill-posed problems is often either limited or not possible for QN methods \cite{smyl2019less,haber2004quasi}, thus effectively limiting applicability in practice. \rev{This is exacerbated in the case of nonlinear problems, where prior information is essential to stabilize reconstructions.}

Motivated by recent advances of data-driven methods for nonlinear inverse problems \cite{martin2017post,hamilton2018deep,hamilton2019beltrami,smyl2020learning,smyl2020optimizing,arridge2020networks,lin2020neural,yoo2019deep,siltanen2020electrical,huuhtanen2020anomaly,agnelli2020classification,candiani2020neural} we propose a new QN approach for updating $\textbf{J}$ based on the singular value decomposition (SVD). First we learn a representative decomposition of the Jacobian from a randomly generated set of training data. We then train a neural network to predict the associated singular values from the model outputs, making it possible to compute an approximate Jacobian without the need for numerical differentiation, resulting in the proposed neural network augmented Quasi-Newton method (NN-QN).
\revtwo{In related studies, researchers have proposed learned SVD frameworks in applications to regularized inverse problems \cite{schwab2019big,schwab2019deep,boink2019learned}} We will first detail the mathematics and framework for the NN-QN approach in \cref{sec:NNaugQN}. We then present an application to  the highly nonlinear inverse problem of electrical impedance tomography in \cref{sec:app2EIT}. \revtwo{In the following,} the approach will be tested with experimental data in \cref{sec:experiments}. Lastly, conclusions will be provided.

\section{A neural network augmented quasi Newton method}\label{sec:NNaugQN}
\subsection{Implementing an approximate Jacobian}
Computing minimizers of $\mathcal{L}(f)$ in \cref{eqn:leastSquaresReg} requires an iterative scheme 
\begin{equation}\label{eqn:updateEqn}
    f_{k+1}=f_k+ \lambda_k \Delta f_k,
\end{equation}
where the search direction $\Delta f_k$  needs to be computed as well as a suitable step-length \rev{$\lambda_k \in (0,2]$   \cite{dennis1996numerical}}. 

Let us first discuss the unconstrained case without a regularization term, i.e. we are solving only the least squares problem $\|g - \mathcal{A}(f)|_2^2$. When using the GN method, we are required to compute $\textbf{J}$ at each iteration in order to attain an unregularized update
$\Delta f_k = \textbf{H}^{-1}\textbf{J}^T (g- \mathcal{A}(f_k)).$
It is worth mentioning that a defining characteristic of the GN method is that higher order terms are ignored when computing the Hessian  $\textbf{H} = \textbf{J}^T\textbf{J}$.
On the other hand, QN methods generally aim to either approximate/update (i) $\textbf{H}$ or $\textbf{H}^{-1}$ using gradient information or (ii) $\textbf{J}$ and use the approximation  of the type $\textbf{H} = \textbf{J}^T\textbf{J}$ \cite{gower2017randomized,martinez2000practical}.
Moreover, when statistical noise information is available, we may include noise weighting in the search direction via a weight matrix $\textbf{W}$.

In this work, we will utilize option (ii) and consider the regularized problem in \cref{eqn:leastSquaresReg}
\rev{which is required for ill-posed problems.} 
Given this choice, we obtain the following weighted description for the update
\begin{equation}
    \Delta f_k = (\textbf{J}^T \textbf{W} \textbf{J} + \Gamma_R )^{-1}(\textbf{J}^T \textbf{W}  (g - \mathcal{A}(f_k)) - \partial R  )
    \label{eqn:QNupdateReg}
\end{equation}
\noindent where $\Gamma_R$ and $\partial R$ correspond to the regularization  and gradient of the regularization term \cite{kaipio2006statistical}. 

In our proposed QN method we utilize the SVD for computing the Jacobian, 
\begin{equation}
    \textbf{J} =\textbf{U}\textbf{S}\textbf{V}^T,
\end{equation}
\noindent with $\textbf{U}$ and $\textbf{V}^T$ containing the left and the right singular vectors, respectively, and $\textbf{S}$ is a diagonal matrix containing the singular values \cite{golub1971singular}.
The crucial assumption for our QN method is that we herein assume that the matrices $\textbf{U}$ and $\textbf{V}^T$ can be reasonably approximated, given a suitable initialization. That is, given an initial $f_0$ we compute the SVD of the corresponding Jacobian i.e. $\textbf{J}(f_0) =\textbf{U}_0\textbf{S}_0\textbf{V}_0^T$. We then fix $\textbf{V}=\textbf{V}_0$ and $\textbf{U}=\textbf{U}_0$ for any approximation of $\textbf{J}$ during the iterative process of \cref{eqn:updateEqn}.
\noindent 
Since $\textbf{S}$ contains the non-negative singular values, it centrally influences the estimation of $\textbf{J}$ and we aim to estimate the updated Jacobian given the fixed singular vectors as
\begin{equation}
    \textbf{J} \approx \textbf{U}_0\textbf{S}\textbf{V}_0^T.
    \label{eqn:Jappx}
\end{equation}
\noindent The challenge is now to obtain a good estimate of the singular values in $\mathbf{S}$ in order to compute the Jacobian (the target quantity) for the QN update. 
To do this, we propose the use of a neural network $\Lambda$ to find an approximate set of singular values $\textbf{S}_\Lambda \approx \textbf{S}$ and obtain the learned estimate for $\textbf{J}$ as
\begin{equation}
    \textbf{J}_\Lambda = \textbf{U}_0\textbf{S}_\Lambda\textbf{V}_0^T.
    \label{eqn:learnedJ}
\end{equation}
\noindent The obtained approximation to the Jacobian can now be readily employed to compute the search direction in \cref{eqn:updateEqn} and forms the basis for the learned QN method proposed here. 
In the next subsection, we will present the neural network (NN) used for predicting and approximating $\textbf{S}_\Lambda$.

\subsection{Learned prediction of the singular values}
The keystone to the viability of the NN augmented Jacobian approximation\revtwo{, as proposed here,} is an accurate and reliable prediction of the singular values in $\textbf{S}_\Lambda$.
We note, that $\textbf{S}_\Lambda$ is a diagonal matrix and hence we are only left with estimating the vector of singular values $\mathrm{diag}(\textbf{S}_\Lambda)$, simplifying the implementation of a learned estimator. 
What remains unclear at present is the selection of the input to the \revtwo{network} 
$\Lambda$ to predict $\mathrm{diag}(\textbf{S}_A)$.
An obvious choice is the model output itself $\mathcal{A}(f)$, since $\textbf{J}$ is a derivative of \rev{$\mathcal{A}(f)$}.
While other choices may be viable, for simplicity we choose this mapping and aim to predict 
\begin{equation}
    \Lambda(\mathcal{A}(f)) = \mathrm{diag}(\textbf{S}_\Lambda)
    \label{Amap}
\end{equation}

\noindent \rev{where $\Lambda$ denotes the functional mapping between NN inputs and predictions.} 
We can then readily form the estimated singular value matrix $S_\Lambda$. The proposed neural network augmented  Quasi-Newton method (NN-QN) can then be obtained by simply substituting \cref{eqn:learnedJ} into the update in \cref{eqn:QNupdateReg}:
 \begin{equation}
 \begin{split}
     \Delta f_k = ((\textbf{U}_0&\textbf{S}_\Lambda\textbf{V}_0^T)^T\textbf{W}(\textbf{U}_0\textbf{S}_\Lambda\textbf{V}_0^T) + \Gamma_R )^{-1} \\
     &((\textbf{U}_0\textbf{S}_A\textbf{V}_0^T)^T\textbf{W}  (g - \mathcal{A}(f_k)) - \partial R  ).
     \label{eqn:QNregLearnedJ}
     \end{split}
 \end{equation}
The major obvious advantage to using the NN-QN method is computational speed, since $\mathcal{A}(f)$ need not be repetitively computed to form the Jacobian after initialization. We also note that Jacobian updates are not susceptible to accumulating round-off errors present in other QN gradient estimators \cite{ding2010investigation}, the NN-QN method \rev{largely} depends on the \revtwo{quality} of $U_0$ and $V_0$ estimates. 
Additionally, it can be expected that $\textbf{S}_\Lambda$ predictions are primarily valid within the space of the training data; but, as we will show in the following, highly general random fields can be successfully used for network training. 
Lastly, we present a summarizing pseudocode of the workflow in Algorithm \ref{alg:NN-QN} as a reference for solving nonlinear inverse problems using the proposed NN-QN method.

\begin{algorithm}
	\caption{NN augmented QN Workflow}
	\label{alg:NN-QN}
	\begin{algorithmic}[1]
	\State Given trained estimator $\Lambda(\mathcal{A}(f)) = \mathrm{diag}(\textbf{S}_\Lambda)$
	\Let{$f_0$}{Initial estimate}
	\State Compute $\textbf{J}(f_0) =\textbf{U}_0\textbf{S}_0\textbf{V}_0^T$
        \Function{NN-QN}{$f_0,\textbf{U}_0,\textbf{V}_0,\Lambda$}
        \Let{$k$}{$0$}
        \While{Stopping criterion $>$ tolerance} 
        \State Compute $\Gamma_R$, $\partial R$, $\mathcal{A}(f_k)$\
        \Let{\rev{$\mathrm{diag}(\textbf{S}_\Lambda)$}}{\rev{$\Lambda(\mathcal{A}(f_k))$}}
        \Let{$\Delta f_k$}{Solve \cref{eqn:QNregLearnedJ}}
        \State Use linesearch to compute $\lambda_k$\
        \Let{$f_{k+1}$}{$f_k+ \lambda_k \Delta f_k$}
        \Let{$k$}{$k+1$}        
        \EndWhile
        \Let{$f_\text{rec}$}{$f_k$}
        \EndFunction{(\Return{$f_\text{rec}$})}
	\end{algorithmic}
\end{algorithm}

\noindent \textbf{Remark:} It is sometimes the case that optimization problems do not require regularization. As such, we would like to note that the proposed NN-QN method also holds for unregularized weighted solutions by using $\Delta f_k = ((\textbf{U}_0 \textbf{S}_\Lambda \textbf{V}_0^T)^T\textbf{W}(\textbf{U}_0\textbf{S}_\Lambda \textbf{V}_0^T) )^{-1}(\textbf{U}_0\textbf{S}_\Lambda\textbf{V}_0^T)^T\textbf{W} (g - \mathcal{A}(f_k))$.  In the case that weighting is not used and assuming $\textbf{J}_\Lambda^T\textbf{J}_\Lambda$ is invertable, we can in fact represent the update in the reduced from 
$\Delta f_k= (\textbf{V}_0 \textbf{S}_\Lambda^{-1}\textbf{U}_0^T) (g - \mathcal{A}(f_k))$.



\section{Neural network augmented Quasi-Newton method for nonlinear inverse problems}\label{sec:app2EIT}

\subsection{Electrical Impedance Tomography}
Herein we demonstrate the proposed NN-QN method using Electrical Impedance Tomography (EIT) \cite{liu2020MMClinear,liu2020Fourier}, a nonlinear inverse problem aiming to estimate the conductivity $\sigma$ provided voltage data $g$. Due to notational convention we write $\sigma$ instead of $f$ for the unknown.
The resultant observation model for EIT is then written as $g = \mathcal{A}(\sigma) + \delta g $ where $\delta g$ is a noise term and $\mathcal{A}(\sigma)$ is the nonlinear forward model characterized by the Complete Electrode Model \cite{somersalo92,vauhkonen99} and solved using finite elements.
In addition to the nonlinear relationship between $\sigma$ and $g$, EIT is ill-posed and susceptible to measurement noise.
As such, regularization is essential to attain suitable reconstructions; to do this, we aim to minimize the following functional characterizing absolute EIT imaging
\begin{equation}
   \mathcal{L}_{\text{EIT}}(\sigma)  =  ||\textbf{L}_e(g - \mathcal{A}(\sigma))||_2^2 + R(\sigma)
\end{equation}
\noindent where $\textbf{L}_e^T \textbf{L}_e = \textbf{W}^{-1}$ and $R(\sigma)$ is the regularization term.
Lastly, we require a linearization point from which to compute the left and right singular \rev{vector} matrices $\textbf{U}_0$ and $\textbf{V}_0^T$.
For this, we simply utilize $\sigma_0 = \sigma_\mathrm{exp}$ where $\sigma_\mathrm{exp}$ is a constant homogeneous estimate based on the expected conductivity, providing a \rev{good estimate and} general starting point.

In the experimental EIT part, data from water tank (non-conductive inclusion localization) and carbon black-modified glass fiber/epoxy laminate (damage localization) tests are considered.
In both cases, 256 measurements are collected using adjacent excitation and measurement following \cite{hauptmann2017open,Tallman2015}.
Geometries for the water tank and laminate tests are circular (radius = 14 cm) and square (width = 9.54 cm), respectively.
These circular and square geometries have 16 electrodes and are discretized using quadratic triangular finite elements meshes consisting of $N_\mathrm{el} = 2034$ and $N_\mathrm{el} = 2528$, respectively.
To attain reliable electrode contact impedances for these finite element models, the best homogeneous estimate was used.
Lastly, in order to demonstrate EIT reconstruction generalizability using the proposed NN-QN method, two different regularization forms are used.
To accomplish this, we use total variation regularization and Laplacian regularization for the water tank and laminate tests, respectively. 

\subsection{Neural network implementation}
The NN architectures used herein employ the same structure, since data sizes for both example problems are identical.
Firstly, we note that the input and output data vectors are of size $\mathbb{R}^{256}$, because the number of singular values is equal to the number of model outputs for the test examples.
Given this, we select a fully connected regression network consisting of three hidden layers composed of 300 neurons each.
The first two layers utilize LeakyReLU activation functions with a negative region slope of $\alpha = 0.1$ while the last hidden layer uses ReLU as activation to ensure non-negative singular values.

\begin{figure}
\centering
\includegraphics[width=3.5cm]{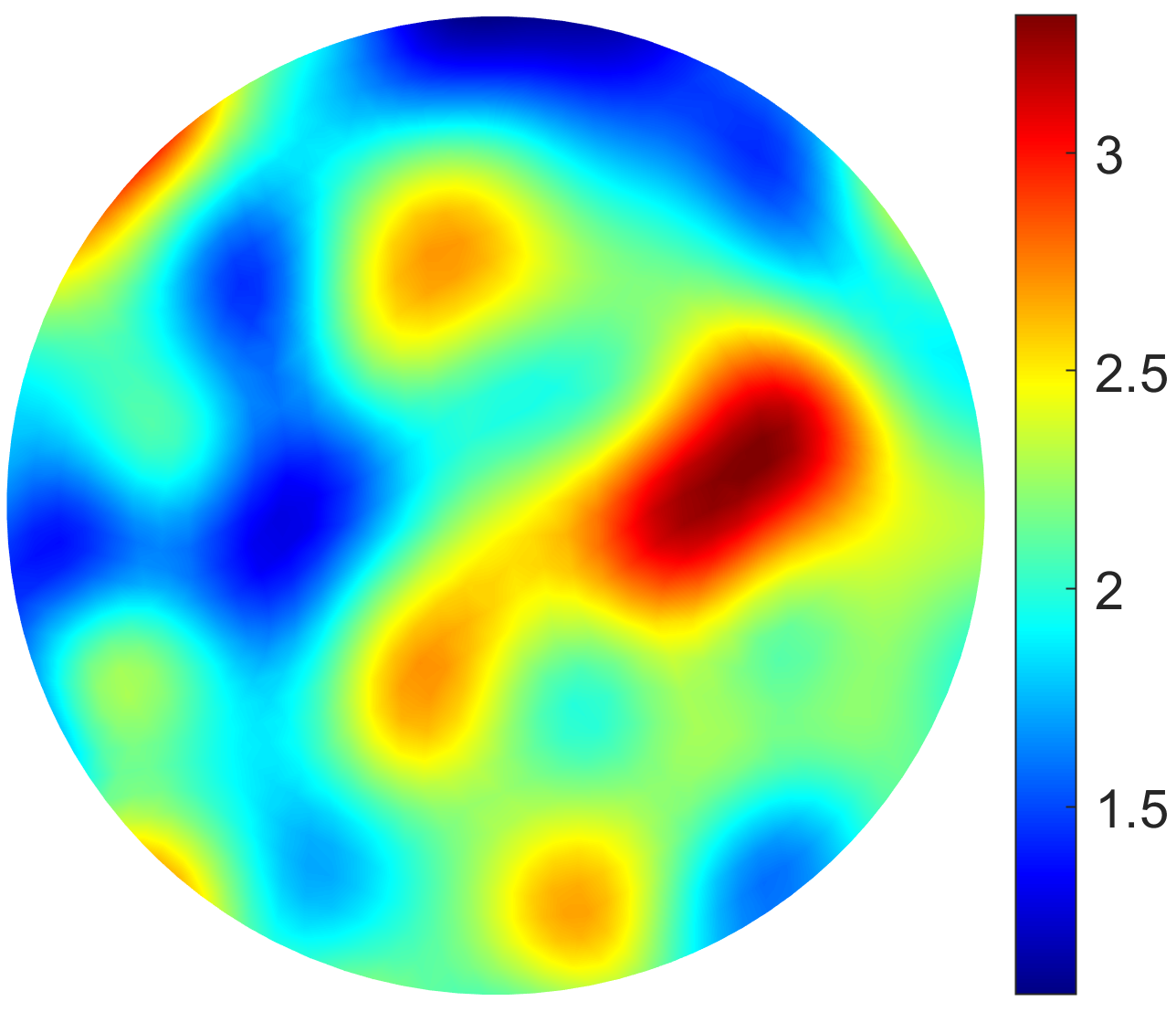} \includegraphics[width=3.75cm]{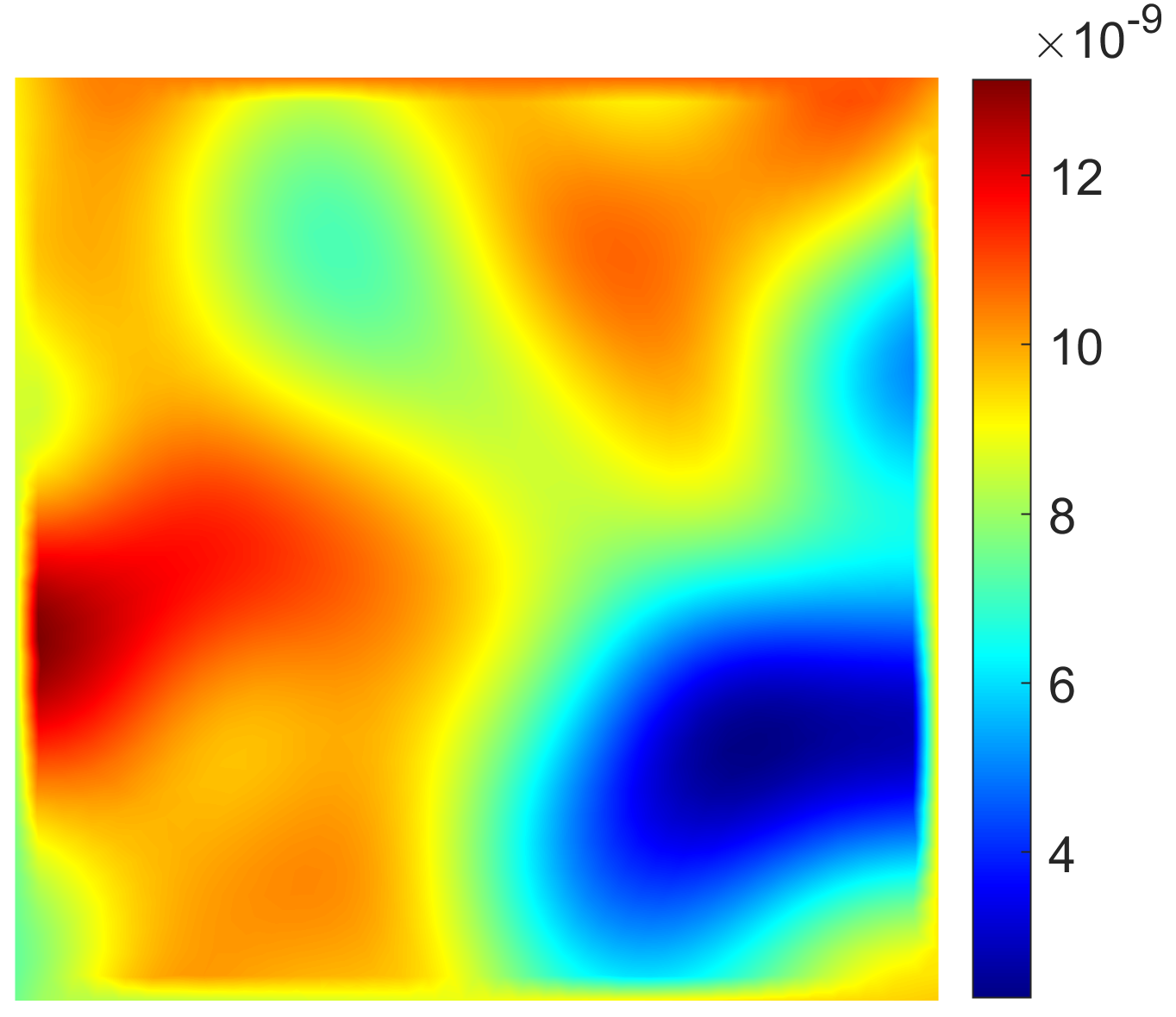}
\caption{Conductivity samples used in generating training data: left, random sample for the water tank geometry (mS/cm); right, random sample used for the laminate geometry (S/m).}
\label{samples}
\end{figure}

To train the networks, $N=5000$ randomized training and $N_t=1000$ validation data sets are used.
The randomized conductivity training samples ${\sigma}_\mathrm{rand}$ were generated using a spatially-weighted Gaussian kernel $\textbf{G}$ perturbed with a random vector $\textbf{r}$, i.e. ${\sigma}_\mathrm{rand} = \textbf{G}^{-1}\textbf{r} + \sigma_\mathrm{exp}$. We note that this is highly general training data, where we ensure the conductivity mean is centered at $\sigma_\mathrm{exp}$ while the range of the random samples varies within two standard deviations of the mean.
To illustrate this, Fig. \ref{samples} shows samples for both EIT geometries.
In optimizing the network using adaptive learning rates (with initial $10^{-5}$), the following loss function is minimized
\begin{equation}
L = \frac{1}{N}\sum_{n=1}^{N} (\mathrm{diag}(\textbf{S}^n) - \mathrm{diag}(\textbf{S}^n_\Lambda))^2 + \kappa \|w\|_2^2,
\label{NNr}
\end{equation}
\noindent where $\mathrm{diag}(\textbf{S}^n)$ and $\mathrm{diag}(\textbf{S}^n_\Lambda)$ represent the $n^\mathrm{th}$ vector of singular values in the training set and output from the NN, $w$ are the network weights, and $\kappa$ is the regularization hyperparameter.
For regularizing the network during training, $\kappa = 0.001$ and a dropout rate between layers of 20\% are used.
To prevent over-fitting, training is terminated when the change in the validation loss was less than $10^{-2}$ between epochs.
\rev{Summarily, learning singular values requires (i) generating model output and singular value samples, (ii) prescribing NN architecture/parameters, and finally (iii) minimizing eq. (\ref{NNr}).}

\section{Experimental results}\label{sec:experiments}
\begin{figure}
\centering
\begin{picture}(220,100)
\put(-0,-8){\includegraphics[width=8.0cm]{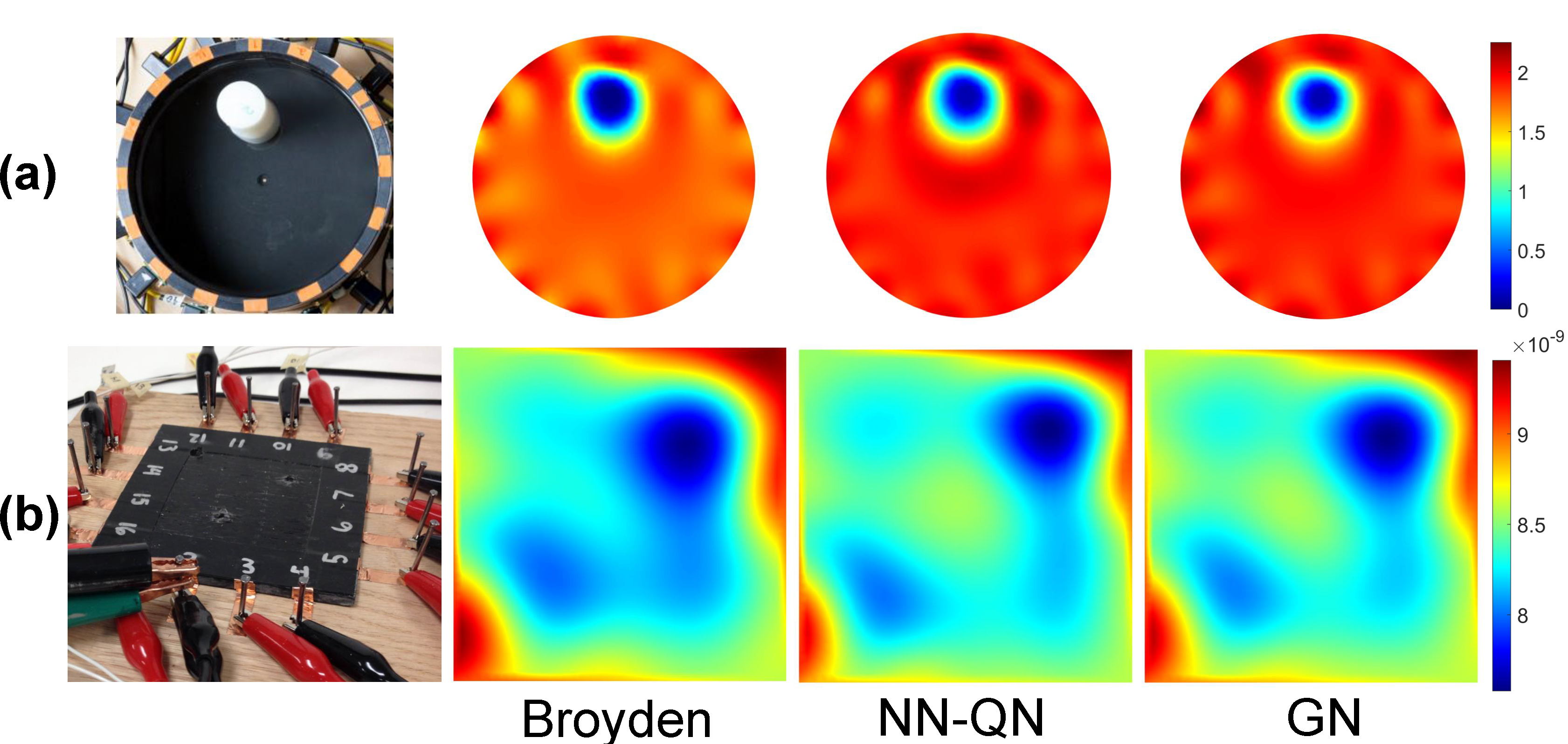}}
\end{picture}
\caption{ EIT reconstructions using (a) water tank data and (b) laminate data. Images in row 1 are experiment photographs while rows 2-4 are computed using the labeled optimization approaches. Colorbar units are mS/cm and S/m for watertank and laminate reconstructions, respectively.}
\label{recons}
\end{figure}

In this section we compare absolute EIT reconstructions from experimental water tank and laminate data as shown in Fig. \ref{recons}.
These reconstructions are computed using (i) the GN method using an accurate perturbed $\textbf{J}$, (ii) the proposed NN-QN method for computing $\textbf{J}$, and (iii) Broyden's classical method for computing $\textbf{J}$.
We note that Broyden's method was selected as it readily supports the integration of $\textbf{W}$ in the inversion.
Additionally, the perturbation method was selected for GN optimization to illustrate computing demands when a semi-analytical algorithm for calculating $\textbf{J}$ is not available (as is often the case).

Water tank reconstructions in Fig. \ref{recons}(a) all show the localization of the non-conductive inclusion well.
\rev{In fact, by visual inspection, NN-QN produces the fewest artifacts near the electrodes but slightly more near the target. We may conclude that NN-QN is comparable to GN, even though GN is based on the more accurate $\textbf{J}$ and $\textbf{H}$.}
Meanwhile, it is also apparent that the background water conductivity estimated using Broyden's method has significant fluctuations near the electrodes and is significantly lower than the accurate GN estimate.
This observation is likely owed to roundoff errors accumulating in successive iterations, consistent with findings in \cite{smyl2019less}.
On the other hand, such roundoff errors do not occur in the proposed NN-QN method, resulting in a more accurate estimation of the background conductivity.

Images reported in Fig. \ref{recons}(b) endeavor to quantitatively localize two through-holes in the composite laminate. These absolute images show that the holes are visible in all cases, although the damage is clearly blurred using Broyden's method.
Due to the small size of the holes relative to the imaged domain and the diffusive nature of EIT, it is clear that the presence of errors/noise in  $\textbf{J}$ and $\textbf{H}$ resulting from the use of Broyden's method significantly impacts the ability of EIT to detect highly localized changes in conductivity.
Comparing the proposed QN-NN and GN reconstructions, the drop in conductivity observed near the holes is visually comparable -- although the NN-QN reconstruction is somewhat smoother relative to the GN reconstruction. 
Note that the conductivity of carbon filler-modified polymers is known to vary appreciably even for a set filler concentration \cite{gungor2015anisotropic}. Hence, conductivity fluctuations away from the holes are expected in Fig. \ref{recons}(b).

A table is provided (Table \ref{comp}), reporting key variables associated with the computational demand of reconstructions presented.
Owing to the fact that the GN method used forward differencing to compute $\textbf{J}$, GN reconstructions took approximately 100-200x longer ($\approx$ 15 minutes) to reach an equal stopping criteria ($||f_{k+1} - f_k||^2_2/||f_{k} ||^2_2 \leq 10^{-2}$) than the Broyden- or NN-based approaches ($\approx$ 5 seconds). We point out that the computation times are after initialization.
Overall, NN-QN and Broyden's method exhibited similar minimization behavior, requiring less time to reach a minimum but requiring more iterations than GN.

\begin{table}
\caption{Computational demand (after Initialization)}
\label{comp}
\centering
\begin{tabular}{@{}lll@{}}
\toprule
Method  & Mean computing time (s) & Mean Iterations (K) \\ \midrule
NN-QN   & 5.22                    & 30        \\
Broyden & 4.93                  & 29.5          \\
GN      & 1858                    & 10.5          \\ \bottomrule
\end{tabular}
\end{table}

\begin{figure}[h]
\centering
\begin{picture}(220,90)
\put(-10,-5){\includegraphics[width=8.5cm]{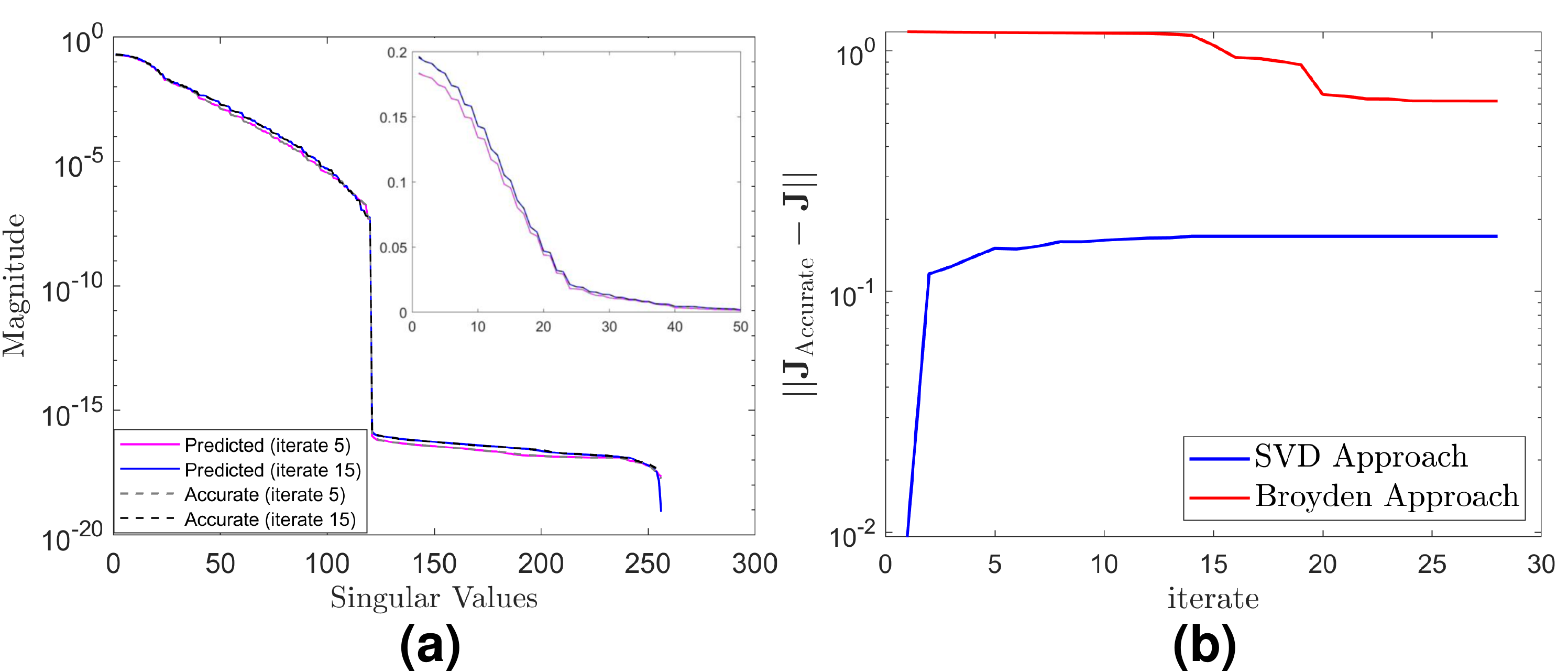}}
\end{picture}
\caption{\rev{Representative verification using EIT water-tank data: (a) Comparison of accurate and predicted singular values {at iterates 5 and 15}; (b) iterative accuracy of the obtained Jacobians compared with Broyden's approach.} }
\label{OptC}
\end{figure}

\rev{We illustrate the effectiveness of singular value predictions, relative to those obtained from an accurate Jacobian $\textbf{J}_\mathrm{Accurate}$, in Fig. \ref{OptC}(a) which compares closely. More importantly, the accurately predicted singular values contribute to a significantly improved approximation of \textbf{J} relative to the Broyden approach, as measured by the differences in matrix norms $||\textbf{J}_\mathrm{Accurate}-\textbf{J}||_2$  shown in Fig. \ref{OptC}(b), which explains the improved reconstruction quality. Specifically, the higher accuracy in early iterates may be beneficial for the final reconstruction quality.
}

\section{Conclusions}
In this work, we presented a new data driven approach to the Quasi-Newton method intended for application to non-linear inverse problems.
The method uses a neural network to map non-linear numerical model outputs to the singular values of the Jacobian matrix, which can be used in forming an approximation of the (otherwise costly) Jacobian.
By doing this, the method falls within a subset of Quasi-Newton methods. 

\rev{Preliminary, the efficacy of the proposed NN-QN method has been demonstrated in the context of an ill-posed non-linear inverse problem - EIT.
It was shown that the use of different priors may be effectively incorporated into the approach.
While these results showed that the NN-QN method may improve reconstructions relative to the use of other QN methods, specific
further research is required to enable a more concrete evaluation.}
\rev{Lastly, it was found that application-specific prior information may be included in the data driven Jacobian estimate via the selection of training data. Whereas very general training data was used herein, more specific distributions may improve solutions to nonlinear inverse problems and is the source of future research.}

\bibliographystyle{unsrt}
\bibliography{bibliography}

\end{document}